\documentclass{article}
\usepackage{amssymb}
\usepackage{amsmath}
\usepackage{amsthm}
\usepackage{graphicx}

\newtheorem{theorem}{Theorem}[section]

\numberwithin{equation}{section}

\begin{document}

\title{Strong solutions to a class of boundary value problems on a mixed Riemannian--Lorentzian metric}
\author{Antonella Marini$^1$ and
Thomas H. Otway$^2$
\\ \\
\textit{$^1$Dipartimento di Matematica, Universit\`{a} di
L'Aquila,}\\ \textit{67100 L'Aquila, Italy } \\ \textit{$^{1,2}$Department of Mathematical Sciences,
Yeshiva University, }\\ \textit{New York, New York 10033}}
\date{}
\maketitle

\begin{abstract}

A first-order elliptic-hyperbolic system in extended projective space is shown to possess strong solutions to a natural class of Guderley--Morawetz--Keldysh problems on a typical domain.

\medskip

\end{abstract}

\section{Introduction}

Although there is a large literature on elliptic--hyperbolic boundary value problems associated with the transition from subsonic to supersonic flow, the literature on boundary value problems that arise from the transition between a Riemannian and a Lorentzian metric is proportionally rather sparse. For reviews, see \cite{OJGP} and \cite{JS}; see also \cite{BC}. In this note we prove the existence of strong solutions to an elliptic--hyperbolic boundary value problem for a lower-order perturbation of the Laplace--Beltrami equation on the extended projective disc $\mathbb{P}^2.$ In this problem the underlying \emph{Beltrami metric} \cite{B, Sw} undergoes a transition from Riemannian to Lorentzian signature along the \emph{absolute}, the curve at infinity, which is the unit circle in $\mathbb{R}^2.$ See, \emph{e.g.}, Sec.\ 9.1 of \cite{H} for discussion. The equations considered here are motivated by at least two topics in spacetime geometry.

First, the Laplace--Beltrami equation on extended $\mathbb{P}^2$ is the hodograph image of the equation for extremal surfaces in Minkowski space $\mathbb{M}^3.$ That equation has the form
\[
\nabla\cdot\left[\frac {\nabla u}{\sqrt{\left \vert1-\vert\nabla u\vert^2\right\vert}}\right] = 0,
\]
where the unknown function $u\left(x,y\right)$ denotes the graph of the surface. See Sec.\ 6.1 of \cite{O4}, and the references cited therein, for discussion.

Second, harmonic fields on extended $\mathbb{P}^2$ have an interpretation as a toy model for waves on certain relativistically rotating cylinders. A rotating, axisymmetric cylindrical solution to the Einstein equations has the general form
\[
ds^2 = -A(r)dt^2+2B(r)d\phi dt+C(r)d\phi^2+D(r)\left(dr^2+dz^2\right),
\]
where $t,z\in\mathbb{R};$ $r\in\mathbb{R}^+;$ and $\phi\in\left[0,2\pi\right].$ This metric has Lorentzian character provided the quantity $AC+B^2$ exceeds zero. Consider the special case of a cylinder rotating at angular velocity $\omega$ and satisfying
\begin{equation}\label{vs}
ds^2=-dt^2+2\omega r^2d\phi dt+r^2\left(1-\omega^2 r^2\right)d\phi^2+F(r)\left(dr^2+dz^2\right),
\end{equation}
where $r<R$ for fixed $R$ and $F$ is a given positive function. The choice $F(r) =\exp\left[-\omega^2 r^2\right]$ yields one of the earliest examples of metrics permitting closed timelike curves; that metric was introduced by van Stockum \cite{VS}. It arises in the context of a rotating, infinitely long cylinder of radius $R$ composed of dust, in which the balance between the centrifugal forces arising from the rotation and the gravitational attraction of the dust provides stability within the surrounding vacuum \cite{Ti}. The possibility of closed timelike curves in this metric occurs when the quantity $r\omega$ exceeds unity. Such curves are associated with causality violation; see, e.g., \cite{LC} for discussion.

Consider the steady case of such a rotation, neglect the height of the cylinder, and make the simplest choice, $F(r)=1.$ The resulting metric has a geometric interpretation as the Beltrami metric on extended $\mathbb{P}^2,$ where $-\omega^2$ emerges as the hyperbolic curvature. While it is much simpler than the metric of eq.\ (\ref{vs}), the Beltrami metric on extended $\mathbb{P}^2$ retains the critical change in signature on the circle $\vert r\omega\vert =1$ (although the dimension reduction flips the Riemannian and Lorentzian regions).

The existence of solutions, having various degrees of regularity, to Laplace--Beltrami equations on metrics having the general form (\ref{vs}) is discussed in Problem 11 of Appendix B to \cite{O4}, in which it is observed that the Laplace--Beltrami operator is not symmetric positive on the stationary case of such a metric. Boundary value problems for weak solutions have however been constructed for that case \cite{OAn}. One expects that a sufficiently strong lower-order perturbation of any symmetric system of equations will be symmetric positive in the sense of Friedrichs \cite{F}. Friedrichs asserted in \cite{F} (but did not prove) that admissible boundary conditions can always been found for a symmetric positive differential equation. But what those admissible conditions are, how to find them, and how to construct a domain on which they will exist, are not known in general and certainly do not follow from Friedrichs' assertion. Such a problem is addressed in this note, for a variant of the Laplace--Beltrami system which is completely integrable and potentially helical. In Theorem 2 of \cite{O1} it is shown that such a solution exists for a particular boundary value problem under stronger hypotheses on the operator than are imposed here, provided the underlying domain has certain abstract properties. Here we construct a class of domains possessing those properties, and construct an explicit class of well posed boundary value problems on those domains. See also \cite{To} for a problem similar to ours, but which is associated to a system of Tricomi type, rather than of Keldysh type as is the case here. In addition, the problem in \cite{To} is posed on a different kind of domain, and with different boundary conditions, than we consider here.

In Sec.\ 2 we formally introduce the system of equations on a domain $\Omega$ and show that the system is symmetric positive on $\Omega,$ under certain hypotheses on the lower-order terms and on the domain. In Sec.\ 3 we impose boundary conditions on $\partial\Omega$ and show that they are admissible in an appropriate sense. We use these properties to prove the existence of strong solutions to an associated boundary value problem by the method of \cite{F}. That existence theorem is stated in Sec.\ 4. Note that the almost-everywhere uniqueness of strong solutions can be shown by what have come to be known as the \emph{Friedrichs inequalities}; see, e.g., Sec. 2.5.1 of \cite{O4} for discussion.

\section{Symmetric positivity of the perturbed operator}

We consider the first-order system
\begin{eqnarray}
\left(1-x^2\right)u_{1x}-xy\left(u_{1y}+u_{2x}\right)+\left(1-y^2\right)u_{2y}\nonumber\\
-\left(2x-\gamma_1u_1\right)-\left(2y-\gamma_2\right)u_2=f_1,\label{LB1}
\end{eqnarray}
\begin{equation}\label{LB2}
u_{1y}-u_{2x}-\left(u_1\Gamma_2-u_2\Gamma_1\right)=f_2.
\end{equation}
Here $U=\left(u_1,u_2\right)$ is an unknown 1-form on a bounded Euclidean domain $\Omega \subset\subset\mathbb{R}^2;$ $\gamma = \left(\gamma_1, \gamma_2\right),$ $\Gamma = \left(\Gamma_1,\Gamma_2\right),$ and $F =\left(f_1,f_2\right)$ are prescribed 1-forms on $\Omega.$ This system is equivalent to a lower-order perturbation of both equations in the Laplace--Beltrami system, on the metric (\ref{vs}) with the choices described in the preceding section, and with the rotational velocity $\omega$ normalized to $1.$ The system (\ref{LB1}, \ref{LB2}) is of elliptic type in the interior of the unit disc $\mathcal{D}_1$ centered at the origin of coordinates in $\mathbb{R}^2.$ It is of hyperbolic type on the $\mathbb{R}^2$-complement of the closure of that disc, and of parabolic type on the boundary of the disc.

For simplicity of exposition, and without compromising the main arguments of the paper, we take $\gamma_2 \equiv 0.$ The extension of the result to nonzero $\gamma_2$ is a matter of algebra. However, we make the important assumption that
\begin{equation}\label{gbound}
\Gamma_1>0\mbox{ and }\vert\gamma_1\vert\geq\frac{1}{4}\frac{\left[xy\Gamma_1+\left(y^2-1\right)\Gamma_2\right]^2}{\left\vert\Gamma_1\left(1-y^2\right)\right\vert}.
\end{equation}
These conditions are necessary and sufficient for the matrix $Q$ in eq.\ (\ref{Qdef}) to be positive on $\Omega.$ Conditions (\ref{gbound}) can be satisfied provided $\Gamma$ is bounded, and the domain is bounded in the $x$-direction, and bounded away from the lines $y^2=1$ in the $y$-direction.  If we take $f_2\equiv 0,$ then eq.\ (\ref{LB2}) is transformed into a condition for complete integrability; see the discussions in Sec.\ 6 of \cite{OJMP} and Secs. 1--3 of \cite{MO}. If $f_2$ is assumed to be nonvanishing, then (\ref{LB2}) becomes a condition for helicity in the sense of \cite{OEd}, Sec.\ 5.4. We assume that $f_1$ is not identically zero, so that we can impose trivial boundary data without obtaining a trivial solution; by linearity an inhomogeneous system with homogeneous boundary data can be shown to be equivalent to a homogeneous differential equation having inhomogeneous boundary data; see Sec.\ 2.6 of \cite{O4}.

Write eqs.\ (\ref{LB1}, \ref{LB2}) as the matrix equation
\begin{eqnarray}
\left(
  \begin{array}{cc}
    1-x^2 & -xy \\
    0 & -1 \\
  \end{array}
\right)
\left(
  \begin{array}{c}
    u_1 \\
    u_2 \\
  \end{array}
\right)_x+
\left(
  \begin{array}{cc}
    -xy & 1-y^2 \\
    1 & 0 \\
  \end{array}
\right)
\left(
  \begin{array}{c}
    u_1 \\
    u_2 \\
  \end{array}
\right)_y+\nonumber\\
\left(
  \begin{array}{cc}
    -2x+\gamma_1 & -2y \\
    -\Gamma_2 & \Gamma_1 \\
  \end{array}
\right)
\left(
  \begin{array}{c}
    u_1 \\
    u_2 \\
  \end{array}
\right)=
\left(
  \begin{array}{c}
    f_1 \\
    f_2 \\
  \end{array}
\right).\label{LBM}
\end{eqnarray}
The system (\ref{LBM}) is not symmetric as a matrix equation, so we solve for $u_{2x}$ in (\ref{LB2}) to obtain
\[
u_{2x}= u_{1y}-u_1\Gamma_2+u_2\Gamma_1.
\]
Substituting this equation into eq.\ (\ref{LB1}) yields
\[
\left(1-x^2\right)u_{1x}-xy\left(2u_{1y}-u_1\Gamma_2+u_2\Gamma_1\right)+\left(1-y^2\right)u_{2y}
\]
\[
-\left(2x-\gamma_1\right)u_1-2yu_2=f_1,
\]
or
\begin{eqnarray}
\left(1-x^2\right)u_{1x}-2xyu_{1y}+\left(1-y^2\right)u_{2y}+\nonumber\\
\left(xy\Gamma_2+\gamma_1-2x\right)u_1-\left(xy\Gamma_1+2y\right)u_2=f_1.\label{LB1a}
\end{eqnarray}
Also, we write in place of (\ref{LB2}),
\begin{equation}\label{LB2a}
\left(1-y^2\right)\left(u_{1y}-u_{2x}-u_1\Gamma_2+u_2\Gamma_1\right)=\left(1-y^2\right)f_2.
\end{equation}
Equation (\ref{LB2a}) is equivalent to (\ref{LB2}) for our choice of $\Omega.$ We obtain the matrix operator $L$ defined by
\begin{eqnarray}
LU= \left(
  \begin{array}{cc}
    1-x^2 & 0 \\
    0 & y^2-1 \\
  \end{array}
\right)
\left(
  \begin{array}{c}
    u_1 \\
    u_2 \\
  \end{array}
\right)_x+
\left(
  \begin{array}{cc}
    -2xy & 1-y^2 \\
    1-y^2 & 0 \\
  \end{array}
\right)
\left(
  \begin{array}{c}
    u_1 \\
    u_2 \\
  \end{array}
\right)_y+\nonumber\\
\left(
  \begin{array}{cc}
    xy\Gamma_2+\gamma_1-2x & -xy\Gamma_1-2y \\
    \Gamma_2\left(y^2-1\right) & \Gamma_1\left(1-y^2\right) \\
  \end{array}
\right)
\left(
  \begin{array}{c}
    u_1 \\
    u_2 \\
  \end{array}
\right)
\equiv A^1U_x+A^2U_y+BU.\label{LBMa}
\end{eqnarray}
Writing
\[
B^\ast = \frac{B+B^T}{2}=\left(
  \begin{array}{cc}
    xy\Gamma_2+\gamma_1-2x & \frac{-xy\Gamma_1-2y+\Gamma_2\left(y^2-1\right)}{2} \\
    \frac{-xy\Gamma_1-2y+\Gamma_2\left(y^2-1\right)}{2} & \Gamma_1\left(1-y^2\right) \\
  \end{array}
\right),
\]
we obtain
\begin{equation} \label{Qdef}
B^\ast-\frac{1}{2}\left(A^1_x+A^2_y\right)\equiv Q,
\end{equation}
where
\[
Q_{11} = xy\Gamma_2+\gamma_1,
\]
\[
Q_{12}=Q_{21}=\frac{-xy\Gamma_1}{2}+\frac{\Gamma_2}{2}\left(y^2-1\right),
\]
and
\[
Q_{22}=\Gamma_1\left(1-y^2\right).
\]
In order for the operator $L$ to be symmetric positive on $\Omega$ in the sense of \cite{F}, we require that $Q_{11}$ be positive, which will be satisfied by condition (\ref{gbound}), and that the matrix determinant
\[
\vert Q \vert = \gamma_1\Gamma_1\left(1-y^2\right)-\frac{1}{4}\left[xy\Gamma_1+\left(y^2-1\right)\Gamma_2\right]^2
\]
also be positive. The latter condition will also follow from (\ref{gbound}) provided $y^2$ is bounded above away from 1 on $\Omega;$ this is insured in the following section.

\section{Admissibility of the boundary conditions}

Define $n_j,$ $j=1,2$ to be the components of the outward-pointing normal. Adopting the summation convention for repeated indices, we write
\[
\beta = n_jA^j =\left(
                  \begin{array}{cc}
                    \left(1-x^2\right)n_1-2xyn_2 & \left(1-y^2\right)n_2 \\
                    \left(1-y^2\right)n_2 & \left(y^2-1\right)n_1 \\
                  \end{array}
                \right).
\]
Writing \cite{O1}
\[
\alpha = \left[-\left(1-y^2\right)\frac{n_2}{n_1}+2xy-\left(1-x^2\right)\frac{n_1}{n_2}\right]n_2,
\]
wherever this object exists we can write $\beta$ in the alternate form
\begin{equation} \label{betalt}
\beta =\left(
                  \begin{array}{cc}
                    -\alpha-\left(1-y^2\right)\frac{n_2^2}{n_1} & \left(1-y^2\right)n_2 \\
                    \left(1-y^2\right)n_2 & \left(y^2-1\right)n_1 \\
                  \end{array}
                \right).
\end{equation}

\begin{figure}
  \includegraphics[width=50mm]{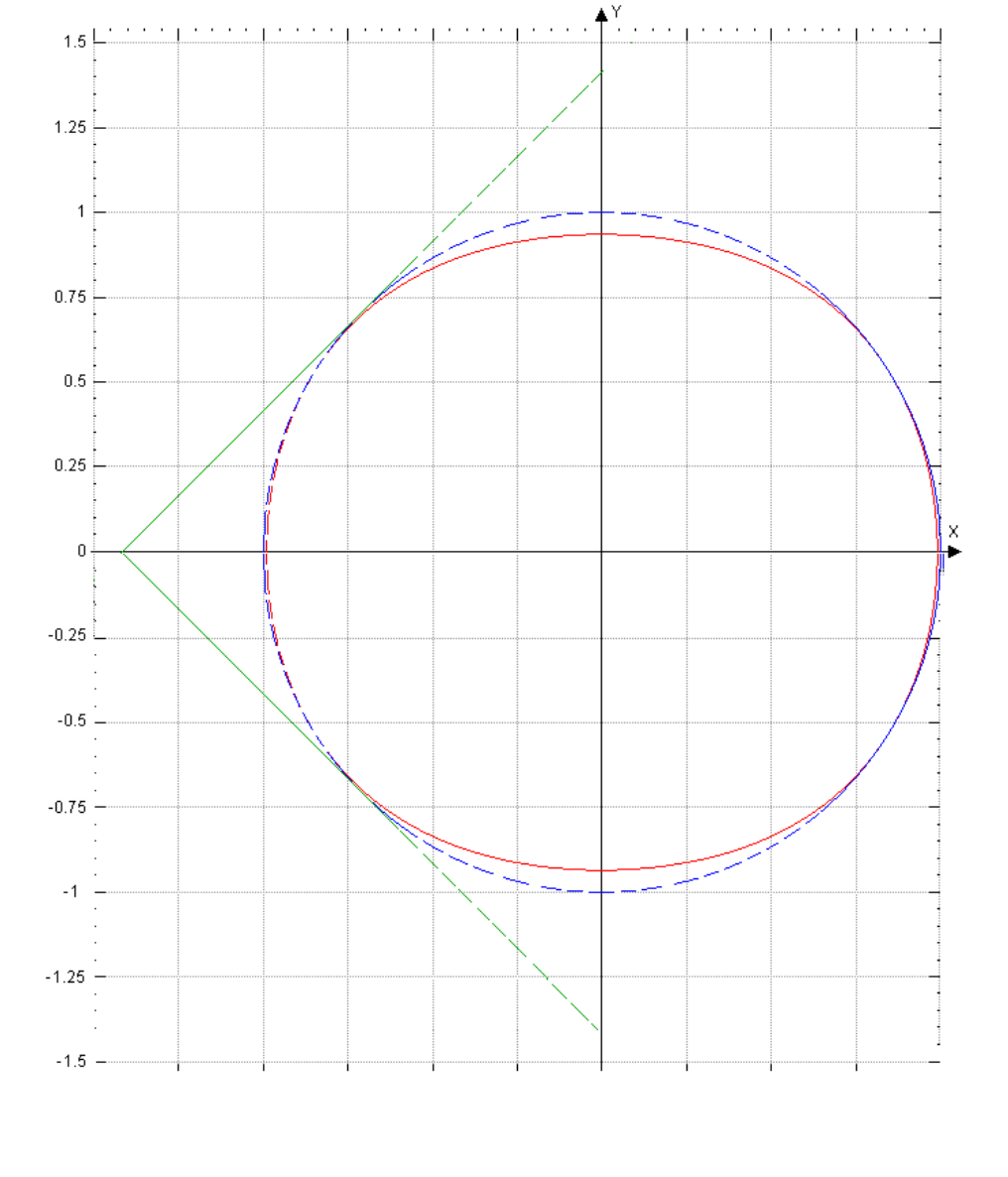}\\
  \caption{Geometry of a typical domain. Hatched lines indicate curves which are not boundary arcs. An arbitrarily small smoothing curve at the corner has been omitted for clarity}  \label{fig1_domain}
\end{figure}

We will consider boundary value problems in the context of the following geometry. Let $\Omega$ include the unit disc $\mathcal{D}_1$ in $\mathbb{R}^2$ centered at the origin of coordinates, truncated at the north and south ``polar caps" by the curves
\begin{equation} \label{c1}
\mathfrak{C}=\left\{\left(x,y\right)\vert y = \pm\left[1-x^2-h\left(x\right)\right]^{1/2}\right\},
\end{equation}
where $h$ is a function chosen so that $0\leq h(x)\leq 1-x^2$ and the graph of $\mathfrak{C}$ is $C^2$ on $\partial{\Omega}.$ The boundary of $\Omega$ is completed in the second and third quadrants by the polar lines $L_1$ and $L_2,$ which are tangent to the unit circle at (for example) the points $\left(-\sqrt{2}/2,\pm\sqrt{2}/2\right).$ Let $h$ vanish at those points. Polar lines have an independent interest in the geometry of extended $\mathbb{P}^2;$ \emph{c.f.} Figure 6.3 of \cite{O4}.  See Figure \ref{fig1_domain}, which illustrates $\overline{\Omega}$ for particular choices of $h,$ $L_1,$ and $L_2.$ Note that the geometry of the domain is rather typical for boundary value problems associated with equations of Keldysh type; \emph{c.f.} Figures 3.2, 3.4, 4.8, and 6.4--6.6 of \cite{O4}.

Make the canonical choice $n_1=dy$ and $n_2=-dx,$ under a ``right-handed" orientation (with the domain interior on the left when the boundary is traversed in the counter-clockwise direction). We find that $\alpha dy=0$ is the equation of characteristic lines to eqs.\ (\ref{LB1a}, \ref{LB2a}) and $\forall\, \left(x,y\right)\in \overline{\mathcal{D}_1},$
\begin{eqnarray}
\alpha dy = -\left(1-y^2\right)dx^2-2xy dxdy-\left(1-x^2\right)dy^2\leq \nonumber \\
-x^2dx^2-2xydxdy-y^2dy^2=-\left(xdx+ydy\right)^2\leq 0.\label{alphan}
\end{eqnarray}

Equation (\ref{alphan}) implies that on any arc $\tau_1,$ contained in $\overline{\mathcal{D}_1},$ on which $dy > 0,$ we have $\alpha\leq 0.$ On $\tau_1$ write $\beta = \beta_+ + \beta_-$ with
\[
\beta_+=\left(
          \begin{array}{cc}
            -\alpha & 0 \\
            0 & 0 \\
          \end{array}
        \right)
\]
and
\[
\beta_-=\left(1-y^2\right)n_2\left(
          \begin{array}{cc}
            -\frac{n_2}{n_1} & 1 \\
            1 & -\frac{n_1}{n_2} \\
          \end{array}
        \right).
\]
Choose the boundary condition $\beta_-U=0,$ that is,
\begin{equation}\label{bct1}
-u_1n_2+u_2n_1=u_1dx+u_2dy=0.
\end{equation}
With this choice of $\beta_+$ and $\beta_-$ we have
\[
\beta_+-\beta_-=\mu=\mu^\ast=\left(
                            \begin{array}{cc}
                              -\alpha+\left(1-y^2\right)\frac{n_2^2}{n_1} & \left(y^2-1\right)n_2 \\
                              \left(y^2-1\right)n_2 & \left(1-y^2\right)n_1 \\
                            \end{array}
                          \right).
\]
Then $\mu_{11}^\ast\geq 0$ and $\vert\mu^\ast\vert=-\alpha\left(1-y^2\right)n_1\geq 0$ for $n_1=dy\geq 0.$

Equation (\ref{alphan}) implies that on any arc $\tau_2,$ contained in $\overline{\mathcal{D}_1},$ on which $dy < 0,$ we have $\alpha\geq 0.$ On $\tau_2$ write $\beta = \beta_+ + \beta_-$ with
\[
\beta_-=\left(
          \begin{array}{cc}
            -\alpha & 0 \\
            0 & 0 \\
          \end{array}
        \right)
\]
and adopt the boundary condition $\beta_-U=0,$ that is,
\begin{equation} \label{bct2}
u_1=0.
\end{equation}
Choose
\[
\beta_+=\left(1-y^2\right)n_2\left(
          \begin{array}{cc}
            -\frac{n_2}{n_1} & 1 \\
            1 & -\frac{n_1}{n_2} \\
          \end{array}
        \right)
\]
so that
\[
\beta_+-\beta_-=\mu=\mu^\ast=\left(
                            \begin{array}{cc}
                              \alpha+\left(y^2-1\right)\frac{n_2^2}{n_1} & \left(1-y^2\right)n_2 \\
                              \left(1-y^2\right)n_2 & \left(y^2-1\right)n_1 \\
                            \end{array}
                          \right).
\]
Then $\mu_{11}^\ast \geq0$ and on any arc within the closure of the unit disc,
\begin{equation} \label{mucirc}
\vert\mu^\ast\vert = \alpha n_1 \left(y^2-1\right)\geq 0.
\end{equation}

On the characteristic lines, $\alpha= 0.$ Then choosing $\beta_\pm$ as on arc $\tau_2,$ $\beta_-$ becomes the zero matrix; so no boundary conditions need to be imposed on the characteristic lines. In this case
\[
\mu^\ast = \left(
             \begin{array}{cc}
               \left(y^2-1\right)\frac{n_2^2}{n_1} & \left(1-y^2\right)n_2 \\
               \left(1-y^2\right)n_2 & \left(y^2-1\right)n_1 \\
             \end{array}
           \right).
\]
Then $\mu_{11}^\ast\geq 0$ as, by construction, $n_1=dy$ is non-positive on the characteristic arcs of $\partial\Omega,$ and $\vert \mu^\ast\vert = 0.$

A boundary value problem in which boundary conditions are imposed everywhere except on the characteristic lines is called a \emph{Guderley--Morawetz} problem. In our problem there is an additional unconventional feature: the ellipticity of the system degenerates on part of the elliptic boundary. Such problems have been studied by Keldysh \cite{K}, and for that reason we refer to the boundary value problem introduced in this section as a \emph{Guderley--Morawetz--Keldysh} problem. However, in distinction to the problem studied in \cite{K}, in this case the elliptic degeneracy plays no role in the analysis. This is an illustration of the powerful type-independence of Friedrichs' method.

\subsection{Singularities and corners}

The matrix $\beta$ in the alternate form (\ref{betalt}) has an apparent singularity at points for which $n_1=0.$ The geometry of the boundary implies that $n_1$ will change sign at no less than two points of $\partial\Omega$ -- and more than two for some choices of $h.$ However, this apparent singularity in (\ref{betalt}) is removable by simply writing out the terms of $\alpha$ and noticing that the singular terms in (\ref{betalt}) cancel additively (for all values of $n_1$).

There is a corner at the intersection of the polar lines (Figure \ref{fig1_domain}). Note that the equations do not change type at this corner, and that the rank of the matrix $\beta$ does not change there. However the conditions of \cite{S1} for regularity at a corner are \emph{not} satisfied; for example, $\vert \mu^\ast\vert$ vanishes at the corner. (See however, \cite{LaP}, an approach which we do not use.) Interpolate an arbitrarily small smoothing curve at the corner.  Then $n_1=dy<0$ at this corner; $n_2=-dx$ is non-negative for $y\geq0$ and non-positive for $y\leq 0.$ Choose
\[
\beta_+ = \left(
                  \begin{array}{cc}
                    \left(1-x^2\right)n_1-2xyn_2 & \left(1-y^2\right)n_2 \\
                    0 & \left(y^2-1\right)n_1 \\
                  \end{array}
                \right)
\]
and $\beta_- = \beta - \beta_+.$ The change of sign in $n_2$ on the $x$-axis is benign, conditions (\ref{bct2}) are satisfied by $\beta_-U,$ and $\mu^\ast\geq 0.$ So the boundary conditions assigned previously to the set $\tau_2$ can be applied to an arbitrarily small smoothing curve at the corner, and such curves are naturally a subset of $\tau_2.$

For all our choices of $\beta_\pm,$ the intersection of the ranges of $\beta_+$ and $\beta_-$ contains only the zero vector. Moreover, the null spaces of $\beta_\pm$ span the restriction to the boundary of the solution space for the system. These properties are required for admissibility in the sense of \cite{F}. In their absence, the boundary conditions are only \emph{semi-admissible}, and only the existence of a weak solution follows from the methods of \cite{F}.

\section{Result}

The arguments of Section 2 and 3 imply that the methods of \cite{F}, which have become standard, can be applied to complete the proof of the following theorem:

\begin{theorem} Let $\overline{\Omega}\subset \mathbb{R}^2$ be the union of the unit disc $\mathcal{D}_1,$ flattened slightly near the poles by the curves $\mathfrak{C}$ given by (\ref{c1}), and the subset of the complement of $\mathcal{D}_1$ which is bounded by the polar lines $L_1$ and $L_2.$ These lines initiate at the points $\left(x_0,y_0\right)$ in the second and third quadrants at which $h\left(x_0\right)=0$ and terminate at an intersection point $\left(x_1,0\right),$ where $x_1<-1.$ (The corner at this intersection can be smoothed to $C^2$ without violating the hypotheses or conclusions of this theorem.) Assume that the prescribed 1-forms $\Gamma_1$ and $\Gamma_2$ do not have blow-up singularities on $\Omega$ and that condition (\ref{gbound}) is satisfied. Then the system (\ref{LB1}, \ref{LB2}), with $\gamma_2\equiv 0$ and $\left(f_1,f_2\right)\in L^2(\Omega),$ supplemented by the boundary condition (\ref{bct1}) on the set $\tau_1\in\overline{\Omega}\backslash\mathfrak{G},$ the boundary condition (\ref{bct2}) on the set $\tau_2\in\overline{\Omega}\backslash\mathfrak{G},$ and no boundary conditions at all on the set $\mathfrak{G}=L_1\cup L_2,$ possesses a strong solution in $\Omega.$
\end{theorem}

\end{document}